\documentclass[10pt]{amsart}

\newtheorem{theorem}{Theorem}

\newtheorem{corollary}[theorem]{Corollary}

\newtheorem{definition}[theorem]{Definition}
\newtheorem{lemma}[theorem]{Lemma}

\newtheorem{proposition}[theorem]{Proposition}
\newtheorem{remark}[theorem]{Remark}

\usepackage{hyperref}

\usepackage[active]{srcltx} 

\def\NN{{\mathbb{N}}}
\def\11{\textbf{$1$}}
\def\CC{{\mathbb{C}}}

\usepackage{enumerate}
\usepackage{amsmath}
\usepackage{amssymb}
\usepackage{soul}

\begin{document}
\begin{center}
    Complete orthogonality preservers between C$^*$-algebras\end{center}
\numberwithin{equation}{section}

\author[J.J. Garc{\'e}s]{Jorge J. Garc\'{e}s}
\address{Departamento de Matem{\'a}tica,
Centro de Ci\^{e}ncias F{\'i}sicas e Matem{\'a}ticas,
Universidade Federal de Santa Catarina, Brazil}
\email{jorge.garces@ufsc.br}

\thanks{Published at Journal of Mathematical Analysis and Applications \url{https://doi.org/https://doi.org/10.1016/j.jmaa.2019.123596}. This manuscript version is made available under the CC-BY-NC-ND 4.0 license \url{https://creativecommons.org/licenses/by-nc-nd/4.0/}}

\title[Complete orthogonality preservers]{}



\date{}
\maketitle

\begin{abstract} We introduce $n$-orthogonality (and completely orthogonality) preserving operators between C$^*$-algebras. Our main theorem states that every completely orthogonality preserving bounded linear mapping between C$^*$-algebras is a weighted TRO homomorphism. We also give several characterisations of TRO homomorphisms among triple homomorphisms. 
 \end{abstract}

\maketitle
\thispagestyle{empty}

\section{Introduction}\label{sec:intro}

 Two elements $a$ and $b$ in a C$^*$-algebra $A$ are said to be \emph{orthogonal}, denoted by $a\perp b$, if $ab^*=b^*a=0.$ A linear mapping $T:A\to B$ between C$^*$-algebras is \emph{orthogonality preserving} if $T(a)\perp T(b),$ whenever $a\perp b,$ $a,b\in A.$
 
 We could trace back the first examples of such mappings to the 
work of E.  Beckenstein,  L.  Narici and  A.R.  Todd on \emph{separating} linear mappings between $C(K)$-spaces (compare with \cite{BeckNarici}). In this context a separating linear mapping (also called \emph{disjointness preserving}) is a mapping that sends functions with disjoint supports to functions with disjoint supports.  W. Arendt proved in \cite{Arendt} that every separating bounded linear mapping between $C(K)$-spaces is a weighted composition operator.  In this setting, the most general description was obtained by K. Jarosz, who was able to describe not necessarily continuous separating linear mappings between $C(K)$-spaces (see \cite{Jarosz}). Jarosz's result was generalised by J.S. Jeang and N.C.
Wong for separating linear mappings between $C_0(L)$-spaces (\cite{JeangWong}).

 In the non-commutative setting several possible generalisations of separating linear mappings arise. Among them are orthogonality preserving and also zero-product preserving linear mappings ($ab=0 \Rightarrow T(a)T(b)=0$).
The first results on zero-product preserving linear mappings on non-abelian C$^*$-algebras are due to L.T. Gardner in \cite{Gardner0}. Zero product preserving linear mappings are termed \emph{disjoint} maps in \cite{Gardner}. Theorem 1 in \cite{Gardner0} together with Theorem 2 in \cite{Gardner} allow to prove that every $2$-positive and disjoint bounded linear mapping between C$^*$-algebras is a \emph{weighted $^*$-homomorphism}. A slightly more general result was obtained by M. Wolff in \cite{Wolff} where it is proved that a symmetric ($T(a^*)=T(a)^*)$ zero-product preserving bounded linear mapping between C$^*$-algebras (being the domain unital) is a \emph{weighted Jordan $^*$-homomorphism}. Wolff's Theorem was later generalised, for zero-product preserving bounded linear mappings, in \cite{CheKeLeeWong}. Zero-product preserving linear mappings coincide with orthogonality preserving linear mappings for abelian C$^*$-algebras and for non-abelian C$^*$-algebras when the mappings are symmetric, however, that is no longer true for non-symmetric maps between non-abelian C$^*$-algebras (see the comments after Theorem 17 in \cite{BurFerGarMarPe}). Zero-product preserving mappings have attracted a lot of attention in the last two decades, resulting on various generalisations and versions in different settings for the just quoted initial results (see for instance \cite{AlaminosZeroP}, \cite{CheKeLeeWong} and \cite{JuChunChinWong}).
  
   Let $A$ be a C$^*$-algebra. We define the triple product of $A$ via $\{a,b,c\}=\frac{1}{2}(ab^*c+cb^*a).$ A linear mapping $T:A\to B$ between C$^*$-algebras is said to be a \emph{triple homomorphism} if $T$ preserves triple products, that is, if $T(\{a,b,c\})=\{T(a),T(b),T(c)\}.$ Triple products play a crucial role in the characterisation of surjective isometries between C$^*$-algebras proved by Kadison in  \cite{Kads}. Indeed, as a consequence of Kadinson's characterisation, it follows that surjective isometries between C$^*$-algebras are triple isomorphisms. Theorem 3.2 in \cite{TripleWong} shows that a bounded linear mapping $T:A\to B$ between C$^*$-algebras is a triple homomorphism if and only if $T$ preserves orthogonality and $T^{**}(1)$ is a partial isometry, thus showing a strong relation between triple products and orthogonality. This relation between triple products and orthogonality becomes even more clear in the general characterisation of orthogonality preserving bounded linear mappings between C$^*$-algebras  obtained by M. Burgos, F.J. Fern{\' a}ndez-Polo,
 J. Mart{\'i}nez Moreno, A.M. Peralta and the author of this note. Indeed, in \cite{BurFerGarMarPe} we prove that every orthogonality preserving bounded linear mapping between C$^*$-algebras is a \emph{weighted triple homomorphism}.

 Orthogonality preserving linear mappings were rediscovered independently from the linear preservers community in \cite{WinterCovering}. Indeed, in \cite[Proposition 4.1.1.]{WinterCovering} completely positive orthogonality preserving linear mappings (there called maps of strict order zero) from a finite dimensional C$^*$-algebra are described. Later in \cite{WinZac}, W. Winter and J. Zacharias  characterise completely positive orthogonality preserving  bounded linear mappings (there termed c.p. order zero maps) between C$^*$-algebras. The main result in \cite{WinZac} shows that any such mapping is a weighted $^*$-homomorphism with positive weight. 
 
 Section 4 in \cite{WinZac} contains some interesting consequences of the main result. Among these is that any amplification of a c.p. order zero map is again an order zero map and, as a consequence, it preserves orthogonality. It is natural to ask whether amplifications of an orthogonality preserving bounded linear mapping also preserve orthogonality. This is not always true. Motivated by that, we shall define $n$-orthogonality preserving bounded linear mappings (also completely orthogonality preserving bounded linear mappings). The main goal of this note is to study this new class of mappings.\smallskip 
 
   A linear mapping $T:A\to B$ between C$^*$-algebras is a TRO homomorphism if $T(ab^*c)=T(a)T(b)^*T(c).$ Clearly, every TRO homomorphism is a triple homomorphism. In Section 2 we give a characterisation of TRO homomorphisms among triple homomorphisms in terms of \emph{right orthogonality}  and \emph{zero-TRO products}. More concretely, we show that the following statements are equivalent:
 \begin{enumerate} [{\rm $1)$}]   
 \item $T$ is a TRO homomorphism,
 \item $T(a)T(b)^*T(c)=0,$ whenever $ab^*c=0$,
  \item $T(a)T(b)^*=0,$ whenever $ab^*=0.$
 \end{enumerate}   
 
Section 3 is devoted to obtain a geometric characterisation of TRO homomorphisms among triple homomorphism. We remind that every TRO homomorphism is completely contractive. Our characterisation shows that a triple homomorphism between C$^*$-algebras is a TRO homomorphism if and only if it is 2-contractive. 

 Inspired by Corollary 4.3 in \cite{WinZac} we introduce in Section 4  $n$-orthogonality and completely orthogonality preserving bounded linear mappings. We are also able to characterise completely orthogonality preserving bounded linear mappings in terms of right orthogonality and zero-TRO products. Our main theorem states that a bounded linear mapping between C$^*$-algebras is completely orthogonality preserving if and only if it is a \emph{weighted TRO homomorphism}.

 In Section 5 we define a (triple) functional calculus for orthogonality preserving bounded linear mappings between C*-algebras. This functional calculus can be regarded as a ternary version of the order zero functional calculus defined in \cite{WinZac}. We also generalise Corollary 4.1 from \cite{WinZac} for contractive completely orthogonality preserving bounded linear mappings.  We apply this generalisation of \cite[Corollary 4.1]{WinZac}  to describe the TRO generated by the range of such an operator.

\section{Preliminaries}

 Throughout the rest of the paper the word ``operator'' will always mean bounded linear mapping. We denote by $B(H)$ the space of operators on a complex Hilbert space $H.$ A TRO is a norm-closed subspace $E$ of some $B(H)$ which is closed for the ternary product $[a,b,c]:=ab^*c.$ Although we shall not define here general JB$^*$-triples, it is worth mentioning that every TRO is a JC$^*$-triple and thus it is also a JB$^*$-triple (compare with \cite{Harris} and \cite{Kaup}). A linear mapping $T:E\to F$ between TROs is a TRO homomorphism (respectively, a TRO anti-homomorphism) if
$T([a,b,c])=[T(a),T(b),T(c)] $ (respectively, $T([a,b,c])=[T(c),T(b),T(a)]$). 

Every C$^*$-algebra is a TRO when endowed with the ternary product $[a,b,c]:=ab^*c$. The (Jordan) triple product of a TRO is defined via $\{a,b,c\}=\frac{1}{2}([a,b,c]+[c,b,a]).$ A norm-closed subspace $I$ of a TRO $E$ is said to be a subTRO if $\{I,I,I\}\subseteq I$ and a triple ideal (respectively, a TRO ideal) if $\{E,I,E\}+\{I,E,E\}\subseteq I$ (respectively, $[E,E,I]+[E,I,E]+[I,E,E]\subseteq I$) holds.  TRO and triple ideals in a TRO coincide, while
TRO, triple and two-sided ideals of a C$^*$-algebras all coincide (see Corollary 5.8 in \cite{Harris}).
     Given an element $a$ in TRO $E$, we define $a^{[1]}:=a,$ $a^{[3]}:=[a,a,a]$ and $a^{[2n+1]}:=[a,a,a^{[2n-1]}]\; (n\in \NN).$ An \emph{odd triple polynomial} on the element $a$ is an element of the form $\sum_{_{k=1}}^{^{n}}\lambda_k a^{[2k-1]},$ for some $n\in \NN$ and $\lambda_k \in \CC.$ We denote by $E_a$ the norm-closed subTRO of $E$ generated by $a$. Clearly, the odd triple polynomials are norm-dense in $E_a.$ It is also known that there exists a (unique) locally compact Hausdorff space $Sp(a)\subseteq (0, \|a\|] $ (called the \emph{triple spectrum} of $a$) such that $Sp(a)\cup \{0\} $ is compact,  $\|a\|\in Sp(a)\cup \{0\} $ and 
there exists a triple isomorphism $\psi : E_a \to C_0(Sp(a))$ mapping $a$ onto the function $t\mapsto t$, where $C_0(Sp(a))$ stands for the Banach space of complex-valued continuous functions on $Sp(h)$ vanishing at $0$ (compare with \cite{Kaup} and take into account that $E$ is a JB$^*$-triple). This allows to define a triple functional calculus on the element $a.$ Given $f\in C_0(Sp(a) )$ we let $f(a):=\psi^{-1}(f)\in E_a$. As an application of the triple functional calculus, there exists a unique element $b\in E_a$ such that $[b,b,b]=a.$ We call this element the \emph{cubic root} of $a$ and denote it by $a^{[\frac{1}{3}]}.$ We can inductively define $a^{[\frac{1}{3^n}]}:=(a^{[\frac{1}{3^{n-1}}]})^{[\frac{1}{3}]}.$ The sequence $(a^{[\frac{1}{3^n}]})$ converges in the weak$^*$-topology of $E^{**}$ to $r(a),$ the range partial isometry of $a$ (see \cite[Lemma 3.3]{EdRut}).

 A linear mapping $T:E\to F$ between TROs is a triple homomorphism if
$T(\{a,b,c\})$ $=\{T(a),T(b),T(c)\}, (a,b,c\in E).$  Since $\{a,b,c\}=\frac{1}{2}([a,b,c]+[c,b,a]),$ it follows that every TRO homomorphism and every  TRO anti-homomorphism is a triple homomorphism. The Jordan product on a C$^*$-algebra is given by $a\circ b =\frac{1}{2}(ab+ba).$  A linear mapping $J:A\to B$ is a \emph{Jordan homomorphism} if $J(a\circ b)=J(a)\circ J(b),$ ($a,b\in A$). Every Jordan $^*$-homomorphism between C$^*$-algebras is a triple homomorphism.

If $A$ is a C$^*$-algebra then, by Sakai's theorem (see Theorem 1.7.8 in \cite{sakai}) the triple product of $A^{**}$ is separately weak$^*$-continuous. This fact together with Goldstine's theorem and weak$^*$-continuity of the second adjoint, allow  to prove that given a triple homomorphism between C$^*$-algebras $S:A\to B$  its second adjoint $S^{**}:A^{**}\to B^{**}$ also is a triple homomorphism. An analogous argument applies to $^*$-homomorphisms, Jordan $^*$-homomorphisms, TRO homomorphisms, etc. Throughout the paper we shall make use of these facts without an explicit mention.

Let $S:M\to N$ be a normal triple homomorphisms between von Neumann algebras and let $r=S(1).$ Then $r$ is a partial isometry and the operator $J=r^*S$ is a normal Jordan $^*$-homomorphism between $M$ and $N$ (see the proof of Corollary 4.3 in \cite{BurFerGarPe}). It also follows that $r^*r$ commutes with $J(M).$  By Stormer's Theorem (see \cite[Corollary 7.4.9]{HaOlSt}), there exist $\phi,\psi:M\to  N,$ such that $ \phi$ is a $^*$-homomorphism, $\psi$ is a $^*$-anti-homomorphism,
  $\phi\perp \psi$ and $J=\psi + \phi.$ It is not hard to see that $rr^*$ also commutes with $\psi$ and $\phi.$
   By using the fact that $rr^*$ also commutes with both $\psi$ and $\phi,$  it can be shown that $r\phi$ is a TRO homomorphism, $r\psi$ is a TRO anti-homomorphism and $r\phi \perp r\psi.$ Finally we have $T=rr^*T=rJ=r\phi+r\psi.$

 \begin{proposition}\label{triplehomDEC} Every normal triple homomorphism between von Neumann algebras decomposes as the orthogonal  sum of a TRO homomorphism and a TRO anti-homomorphism. $\hfill \Box$\end{proposition}

A celebrated result by Bre\v sar (see Theorem 3.14 in \cite{CheKeLeeWong}) states that a Jordan isomorphism between C$^*$-algebras that preserves zero products is an algebra isomorphism. In \cite{PeraltaLocalRep} A.M. Peralta   generalises Bre\v sar's theorem by showing that every bounded Jordan homomorphism between C$^*$-algebras that preserves zero products is a homomorphism. We shall next obtain a similar characterisation of TRO homomorphisms among triple homomorphisms. We remind that Corollary 18 in \cite{BurFerGarMarPe} states that an operator between C$^*$-algebras preserves orthogonality if and only if it preserves zero-triple products (that is, $\{T(a),T(b),T(c)\}=0,$ whenever $\{a,b,c\}=0$).  

  \begin{definition} Let $T:A \to B$ be a linear mapping between C$^*$-algebras. We say that $T$ preserves zero-TRO products if $[T(a),T(b),T(c)]=0$ whenever $[a,b,c]=0.$ \end{definition}

It is clear that every TRO homomorphism preserves zero-TRO products. Let $S:M_2(\CC)\to M_2(\CC),\;S(a)=a^t.$ Then $S$ preserves orthogonality ($S$ is a triple homomorphism),
  however, it does not preserve zero-TRO products. Indeed, take $a=\left [\begin{array}{ll}
         1 & 0\\
        0 & 0 \end{array} \right ]$ and $b=\left [\begin{array}{ll}
         0 & 0\\
        1 & 0 \end{array} \right ].$ Then $[a,a,b]=0$ but $[S(a), S(a),S(b)]=\left [\begin{array}{ll}
         0 & 1\\
        0 & 0 \end{array} \right ].$\medskip
        
         As Lemma 1 in \cite{BurFerGarMarPe} shows, two elements $a,b$ in a C$^*$-algebra are orthogonal  if and only if
 $\{a,b,b\}=0$. We show a similar useful characterisation for \emph{right orthogonality}. 
 
  \begin{lemma}\label{lemma TRO-zero elements} Let $a,b$ be two elements in a C$^*$-algebra $A.$ Then $ab^*=0$ if and only if $[a,b,b]=0.$
  \end{lemma}

  \begin{proof}
   $\Rightarrow)$ It is straightforward.
    $(\Leftarrow $ Let us fix $a,b$ in $A$ such that $ab^*=0$ and let $C^*(b^*b)$ be the C$^*$-subalgebra of $A$ generated by $b^*b.$ We easily check that  $a(b^*b)^n=0$ for all $n\in \NN.$ It is easy to see now that $ap(b^*b)=0$ for any polynomial with  $p(0)=0$. By norm density of these polynomials in  $C^*(b^*b)$ we have $ax=0$ for every $x\in C^*(b^*b).$ In particular $a |b|^{\frac{1}{n}}=0$  holds for all $n\in \NN$. 
    By separate weak$^*$-continuity of the product in $A^{**}$ we have $a s(|b|)=0,$ where $s(|b|)$ is the support projection of $|b|.$ Now observe that $s_r(b)=s(|b|),$ where $s_r(b)$ is the right support of $b.$ Finally we have $$ab^*=as_r(b)b^*=as(|b|)b^*=0. $$
    \end{proof}
    
    Right (left) orthogonality preserving operators have also been studied in the contexts of C$^*$-algebras and preduals of von Neumann algebras (see \cite{LauWong}). The reader can check the survey \cite{JuChunChinWong} and references therein for the main results on zero-product, orthogonality or right (left) orthogonality preserving operators between C$^*$-algebras and some other related structures. In the next result we show that TRO-zero products and right orthogonality allow to characterise TRO homomorphisms among triple homomorphisms. Although this might be known by experts, we include here a short proof.

  \begin{proposition}\label{char TRO-hom} Let $S:A\to B$ be a triple homomorphism between C$^*$-algebras. The following statements are equivalent:
  \begin{enumerate}[$1)$]
  \item $S$ is a TRO homomorphism,
  \item $S$ preserves TRO-zero products,
  \item $S(a)S(b)^*=0,$ whenever $ab^*=0.$
   \end{enumerate}\end{proposition}

   \begin{proof} $1)\Rightarrow 2)$ Is clear and   $2)\Rightarrow 3)$ follows from Lemma \ref{lemma TRO-zero elements}.\medskip

    $3)\Rightarrow 1)$ The proof of Theorem 3.6 in \cite{JuChunChinWong} can easily be adapted for right orthogonality to prove that $S$ satisfies $$S(ab^*)S^{**}(1)^*=S(a)S(b)^*\; ( a,b \in A). $$
    
     Since $S$ is a triple homomorphism then $r:=S^{**}(1)$ is a partial isometry and $S(a)=rr^*S(a)=S(a)r^*r$. We have
     $$ S(ab^*c)=S(ab^*c)r^*r=S(ab^*)S(c^*)^*r=S(ab^*)r^*rS(c^*)^*r$$
    $$=S(a)T(b)^*rS(c^*)^*r=S(a)S(b)^*\{S^{**}(1),S(c^*),S^{**}(1)\}$$ $$=S(a)S(b)^*S(\{1,c^*,1\})=S(a)S(b)^*S(c),$$ which proves that $S$ is a TRO homomorphism.
    
\end{proof}
   
   \section{A geometric characterisation of TRO homomorphisms}

   Let $T:A\to B$ be a linear mapping between C$^*$-algebras. The $n$-th amplification of $T,$ denoted by $T_n,$ is the linear mapping $T_n: M_n(A)\to M_n(B)$ given by $T_n((a_{i,j}))=(T(a_{i,j})).$ $T$  is said to be an $n$-isometry if $T_n$ is an isometry. $T$ is a
\emph{complete isometry} if $T$ is an $n$-isometry, for all $n\in \NN.$ It is well known that surjective isometries between C$^*$-algebras are precisely triple isomorphism while complete surjective isometries are TRO isomorphisms (see \cite{Kads} and \cite[Proposition 2.1]{HAM}). $T$ is said to be $n$-contractive if $T_n$ is contractive and \emph{completely contractive} if each $T_n$ is contractive. While triple homomorphisms between C$^*$-algebras are always contractive they are not in general completely contractive (for instance, the transpose on matrix algebras is not completely contractive as proved in \cite{tomiyama}). On the other hand, every TRO homomorphism is completely contractive (\cite[Proposition 2.1]{HAM}). A natural question is whether every completely contractive triple homomorphism must be a TRO homomorphism. We shall devote the current section to answer this question.

  Let $S:A\to B$ be a triple homomorphism between C$^*$-algebras with $A$ abelian. Then $S^{**}:A^{**}\to B^{**}$ also is a triple homomorphism. By Proposition \ref{triplehomDEC}, $S^{**}=\phi+\psi,$ where $\phi:A^{**}\to B^{**}$ is a TRO homomorphism, $\psi:A^{**}\to B^{**} $ is a TRO anti-homomorphism and $\phi(A^{**})\perp \psi(A^{**}).$ Let us check that $\psi$  is a TRO homomorphism. For $a,b,c\in A^{**},$ we have $\psi(a)\psi(b)^*\psi(c)=\psi(cb^*a)=\psi(ab^*c) ,$  and thus $\psi(ab^*c)=\psi(a)\psi(b)^*\psi(c).$ Hence $S^{**}$ is an orthogonal sum of two TRO homomorphisms whence $S^{**}$ is a TRO homomorphism. As a consequence, $S$ is a TRO homomorphism.

\begin{lemma}\label{l from abelian is TRO} Every triple homomorphism from an abelian C$^*$-algebra is a TRO homomorphism.$\hfill \Box$
\end{lemma} 
 
A C$^*$-algebra $A$ is self-opposite if $A$ is $^*$-isomorphic to $A^{op}$ (equivalently,  if there exists a $^*$-anti-automorphism on $A$).  Examples of self-opposite C$^*$-algebras are $K(H)$ and $B(H)$ for $H$ finite or infinite dimensional. More generally every groupoid C$^*$-algebra is self-opposite (compare with \cite{bussim}).  Examples of C$^*$-algebras which are not self-opposite can be found in \cite{connes}.

 Let $E$ be a TRO. We denote by $E^{op}$ the TRO with underlying Banach space $E$ and  product given by $[a,b,c]_{op}:=[c,b,a].$ A TRO is said to be self-opposite if there exists a TRO isomorphism between $E$ and $E^{op}$  (equivalently, if there exists a TRO anti-isomorphism on $E$). Let $A$ be a C$^*$-algebra realised as a TRO. If $A$ is self-opposite (as a TRO) then there exists a TRO isomorphism $S:A\to A^{op}.$  By Corollary 2.2 in \cite{HAM} there exists a $^*$-isomorphism from $A$ onto $A^{op}.$ 

  \begin{lemma} \label{l techselfopp1}
A C$^*$-algebra is self-opposite if and only if it is self-opposite as a TRO.

Let $E,F$ be TROs and $S:E\to F$ be injective. Then we have:
\begin{enumerate}[$1)$]

\item If $E$ is self-opposite  and $S$ is TRO homomorphism  then $S(E)$ is a self-opposite TRO,

\item If $E$ is self-opposite and  $S$ is TRO anti-homomorphism  then $S(E)$ is a self-opposite TRO,

\item If $S$ is a TRO anti-isomorphism and $E$ or $F$ is self-opposite, then there exists a TRO isomorphism between $E$ and $F$.
\end{enumerate}
\end{lemma}

\begin{proof} 
$1)$ 
Let $S: E\to F$ be an injective TRO homomorphism into another TRO $F.$ Then $S(E)$ is a TRO ($S$ has closed range by Corollary 7 in \cite{KapsTriples}). As a consequence, $S:E\to S(E)$ is a TRO isomorphism. It is easy to realise that
$\overline{S}:E^{op}\to S(E)^{op}, a\mapsto \overline{S}(a)=S(a)$ also is a TRO isomorphism. Now let $\widetilde{S}: S(E)\to E, \widetilde{S}=S^{-1}_{|S(E)}.$ Clearly, $\widetilde{S}$ is a TRO isomorphism. Since $E$ is self-opposite, there exists a TRO isomorphism $T:E\to E^{op}.$ Finally, $\overline{S}\circ T\circ \widetilde{S}:S(E)\to S(E)^{op}$ is a TRO isomorphism.

$2)$ Now suppose that $S: E\to F$ is an injective TRO anti-homomorphism into another TRO $F.$ Then $S':E\to F^{op}, S'(a)=S(a)$ is a TRO homomorphism. By $1)$ the TRO $S'(E)=S(E)$, when endowed with the product of of $F^{op},$ is self-opposite. Equivalently, $S(E)^{op}$ is self-opposite (when endowed with the product of $F$). Thus $S(E)$ is self-opposite.

$3)$  We may assume that $F$ is self-opposite. Let $\phi:F\to F$ be a TRO anti-isomorphism. Then $\phi \circ S:E\to F$ is a TRO isomorphism.
\end{proof}

A linear map between C$^*$-algebras $T:A\to B$ is said to be \emph{completely
bounded} with completely bounded norm given by $$\|T\|_{cb}:=\sup \{ \|T_n\|: n\in \NN \} $$ if the latter is finite.

 Let $\theta(m):M_m(\CC)\to M_m(\CC)$ be the transpose map. Then by \cite[Theorem 1.2]{tomiyama} \begin{equation}\label{label transposefinite} \|\theta(m)_{_n}\|= \left\{
\begin{array}{c l}
 n, & \mbox{ if } n\leq m\\
 m, & \mbox{ if } n>m
\end{array}
\right.  \end{equation} We see that $\theta(m)$ is completely bounded with $\|\theta(m)\|_{cb}=m.$ If $H$ is an infinite dimensional Hilbert space, we write $\theta$ for the transpose on $B(H).$ In this case $\| \theta_n\|=n$ (compare with \cite{tomiyama}). Thus $\theta$ is not completely bounded.\smallskip

 Let $S:B(H)\to B(H)$ be a triple isomorphism (equivalently, a surjective linear isometry). Then $S$ is either of the form $S(a)=v^*au$ (a TRO isomorphism) or of the form 
 $S(a)=v^*a^tu$ (a TRO anti-isomorphism) where $u,v$ are two unitaries in $B(H)$ (compare with \cite{franzoni}).
 In the first case $S$ is a complete surjective isometry with $\| S_n\|=1$ for all $n\in \NN.$ Suppose that  $S(a)=v^*a^t u.$ Since the mapping $a\mapsto v^*au$ is a complete isometry, we have $\|S_n(a_{i,j})\|=\|v^* (a_{i,j})^t u\|=\|  (a_{i,j})^t \|,$ $(a_{i,j})\in M_n(B(H)).$

\begin{proposition}\label{p CCfactors} Let $S:B(H)\to B(H)$ be a TRO anti-isomorphism. We have:
\begin{enumerate}[$1)$]
\item If $H$ is finite dimensional then $\| S_n\|=\|\theta (m)_n \|,$ where $m=dim(H),$

\item If $H$ is infinite dimensional then $\| S_n\|=\|\theta_n \|=n.$
\end{enumerate}
$\hfill \Box$
\end{proposition}

 By Lemma \ref{l from abelian is TRO}  we know that every triple homomorphism from an abelian C$^*$-algebra is a TRO homomorphism, and hence it is completely contractive. We shall see that this is never true for non-abelian C$^*$-algebras.
 
\begin{theorem} \label{t cctriplehom} Let $S:M\to N$ be a normal TRO anti-homomorphism between von Neumann algebras. Suppose that $M$ is not abelian. Then $\|S_n\|\geq 2.$ 

Moreover, if $M$ contains a copy of $M_m(\CC)$ then $\|S_n\|\geq \|\theta (m)_n\| .$ 
\end{theorem}
\begin{proof} There is no loss of generality in assuming that $S$ is injective. Indeed, $I=ker(S)$  is a weak$^*$-closed ideal of $M$ and  $M=I\oplus^{\infty}I^{\perp}.$ It is not hard to see that $M_n(M)=M_n(I)\oplus^{\infty}M_n(I^{\perp}).$ We can further assume that $I^{\perp}$ is not abelian for, if it was abelian then $S$ would be  a TRO homomorphism. Thus $\widetilde{S}=S_{|I^{\perp}}$ is an injective TRO anti-homomorphism with 
  $\| \widetilde{S}_n\|=\| S_n\|.$    
  
  We assume now that $S$ is injective. Let $p,q$ be two non-commuting projections in $M.$ Then the von Neumann algebra generated by $p$ and $q,$ W$^*(p,q),$ is of the form $\CC^n\oplus C(K)\otimes M_2(\CC)$ with $n\leq 4$ and $K$ Stonean (compare with Theorem 1.41 in \cite[Chapter 8]{takesaki})). Since W$^*(p,q)$ is not abelian, then it has non-trivial non-abelian part. In particular, W$^*(p,q)$ contains a copy of $M_2(\CC).$ As a consequence, $M$ contains a $^*$-isomorphic copy of $M_2(\CC).$ Let $\widetilde{S}:=S_{|M_2(\CC)}:M_2(\CC)\to N.$  Since $S$ is injective, then $\widetilde{S}\neq 0.$
 By Lemma \ref{l techselfopp1}, $\widetilde{S}(M_2(\CC))$ is a self-opposite TRO which is TRO isomorphic to $M_2(\CC).$ There exists a TRO isomorphism $\phi:S(M_2(\CC))\to M_2(\CC).$ Let us define $S':M_2(\CC)\to M_2(\CC)$ by $S'=\phi \circ \widetilde{S}.$ Then $S'$ is a TRO anti-isomorphism. Clearly, $S_n'=\phi_n \circ \widetilde{S}_n.$ Since $\phi$ is a complete isometry, we have $$\|S'_n(a)\|=\|\phi_n(\widetilde{S}_n(a))\|=\|\widetilde{S}_n(a)\|. $$ Since $S'$ is a TRO anti-isomorphism on $M_2(\CC)$, we have $\|S'_n\|=\|\theta(2)_n\|=2$ (see  Proposition \ref{p CCfactors}). Finally we have
$$\|S_n\| \geq \|(S_{|M_2(\CC)})_n \| =\| \widetilde{S}_n\| =2.$$

 If $M$ contains a copy of $M_m(\CC),$ then analogous computations allow us to prove that for $\widetilde{S}:=S_{|M_m(\CC)}$ we have 
$\|S_n\| \geq \| \widetilde{S}_n\| =\|\theta (m)_n\| .$

 \end{proof}

Let $S:M \to N$ be a normal triple homomorphism between von Neumann algebras. By Proposition \ref{triplehomDEC} there exists a TRO homomorphism  $\phi$ and a TRO anti-homomorphism $\psi$ such that $S=\phi + \psi$ and $\phi(M)\perp \psi(M)$. It is not hard to see that $S_n=\phi_n + \psi_n$ with $\phi_n(x)\perp \psi_n(y),\; x,y \in M_n(M). $ We have
$$\|S_n(x)\|=\max \{\| \phi_n(x)\|, \|\psi_n(x)\| \}, x\in M_n(M).$$  If $S$ is not a TRO homomorphism, that is if $\psi\neq 0$, then by Theorem \ref{t cctriplehom}, $\|S_n\|=\|\psi_n(x)\|\geq 2.$ Now let $S:A\to B$ be a  triple homomorphism between C$^*$-algebras. Then $S^{**}:A^{**}\to B^{**}$ is a  normal triple homomorphism and $\|S_n\|=\|S_n^{**}\|\geq 2.$ Thus if $S$ is a  triple homomorphism and $\|S_n\|<2$ for some $n\in \NN,$ then $S$ is a TRO homomorphism.

\begin{corollary}\label{c ccTriple cstar}
Every 2-contractive triple homomorphism between C$^*$-algebras is a TRO homomorphism.$\hfill \Box$
\end{corollary}

\section{Completely orthogonality preserving operators}

Let $T:A\to B$ be a completely positive (in short c.p.) linear mapping between C$^*$-algebras. Following \cite{WinZac} we say that $T$ has \emph{order zero} if $T(a)T(b)=0,$ whenever $ab=0,$ with $a,b\in A_{+}.$ Since $T$ is, in particular, a positive map, that is equivalent to saying that $T$ preserves orthogonality on positive elements. The main theorem of \cite{WinZac} states that every completely positive order zero map between C$^*$-algebras is a weighted $^*$-homomorphism.
 Among the consequences of this characterisation of c.p. order zero maps, is that every amplification of a c.p. order zero map is an order zero map. In particular, amplifications of  order zero maps preserve orthogonality. We notice that amplifications of an orthogonality preserving operator may not preserve orthogonality. \smallskip
 
 Let $\theta(2)$ be the transpose on $M_2(\CC).$ Since $\theta(2)$ is a triple homomorphism then it preserves orthogonality. However $\theta(2)_2$ is not orthogonality preserving. Indeed, let us take $x=\left [\begin{array}{llll}
         1 & 0&1&0\\
        0 & 0&0&0\\
        1 & 0&1&0\\
        0 & 0&0&0\end{array} \right ],$   $y=\left [\begin{array}{llll}
        0 & 0&0&0\\
        1 & 2&-1&0\\
        0 & 0&0&0\\
        2 & 2&-2&0\end{array} \right ] \in M_2(( M_2(\CC)).$
        
         Then $x\perp y,$ however 
        $\theta(2)_2(x)$ and $\theta(2)_2(y)$ are not orthogonal elements since $$(\theta(2)_{2}(x))^*\theta(2)_{2}(y)= \left [\begin{array}{llll}
         0 & 3&0&-3\\
        0 & 0&0&0\\
        0 & 3&0&-3\\
        0 & 0&0&0\end{array} \right ] .$$

  Our first step is to study when  amplifications of a Jordan homomorphism (respectively, a triple homomorphism) are again a Jordan homomorphism (respectively, a triple homomorphism). This is actually well known. We state the results as a lemma for completeness reasons.
 
\begin{lemma}\label{l 2-jordantriple} Let $S:A\to B$ be a linear mapping between C*-algebras. Then we have

\begin{enumerate}[$1)$]

\item $S_2$ is a Jordan homomorphism if and only if $S$ is multiplicative.

\item $S_2$ is a  triple homomorphism if and only if $S$ is a TRO homomorphism.

\item If $S$ is a TRO homomorphism then $S_n$ also is a TRO homomorphism, for every $n\in \NN.$
\end{enumerate}
\end{lemma} 
\begin{proof}
See Lemma 2.3 in \cite{HAM} for the statement regarding triple homomorphisms. The one regarding Jordan homomorphisms follows by similar arguments.
\end{proof}

 In the rest of the section we shall be concerned with studying those orthogonality preserving operators whose amplifications also preserve orthogonality. 

 \begin{definition} We say that $T:A\to B$ is \emph{$n$-orthogonality preserving} if $T_n$ is orthogonality preserving. $T$ is said to be
 \emph{completely orthogonality preserving} (in short \emph{c.o.p.}) if  $T_n$ is orthogonality preserving for every $n\in \NN$.\end{definition}

  Every TRO homomorphism is completely orthogonality preserving. Let $S:A \to B$ be a 2-orthogonality preserving triple homomorphism.  Since $S$ is a triple homomorphism, then $r=S^{**}(1)$ is a partial isometry. It is easy to see that
   $S_2^{**}(\left [\begin{array}{ll}
         1 & 0\\
        0 &1  \end{array} \right ])= \left [\begin{array}{ll}
         r & 0\\
        0 &r  \end{array} \right ]$ also is a partial isometry. Thus $S_2$ preserves orthogonality and $S_2 (1_{_{M_2(A^{**})}} )$ is a partial isometry. By \cite[Theorem 3.2]{TripleWong} (see also \cite{BurFerGarMarPe} or \cite{BurFerGarPe}) $S_2$ is a triple homomorphism. By Lemma \ref{l 2-jordantriple},  $S$ is a TRO homomorphism.\smallskip
        
        Moreover, if $S$ is a Jordan $^*$-homomorphism, then $r$ is a projection in $B^{**}$ and $\left [\begin{array}{ll}
         r & 0\\
        0 &r  \end{array} \right ]$ is a projection in $M_2(B)^{**}.$ By Lemma 3.3 in \cite{Wolff} $S_2$ is a Jordan $^*$-homomorphism and hence $S$ is a $^*$-homomorphism by Lemma \ref{l 2-jordantriple}.

\begin{proposition} \label{p triple2OP} Let $S:A\to B$ be a a triple homomorphism  between C$^*$-algebras. Then we have:
 
 \begin{enumerate}[$1)$]
 \item $S$ is a TRO homomorphism if and only if $S$ is 2-orthogonality preserving.
 
 \item $S$ is a $^*$-homomorphism if and only if $S$ is a 2-orthogonality preserving Jordan $^*$-homomorphism.
 \end{enumerate}
 Moreover, every 2-orthogonality preserving triple homomorphism is completely orthogonality preserving.
 $ \hfill \Box$
  \end{proposition}

The characterisation of c.p. order zero maps obtained by W. Winter and J. Zacharias relied on Wolff's description of symmetric orthogonality preserving operators and Stinespring Theorem. It is worth mentioning that this characterisation could also be derived from Gardner's results (compare with \cite{Gardner}). Our characterisation of c.o.p. operators  relies on the characterisation of orthogonality preserving operators obtained in \cite{BurFerGarMarPe} which we now state for later use.

 \begin{theorem}\label{t OPcstar} Let $T: A\to B$ be an orthogonality
 preserving operator between two C*-algebras. For $h= T^{**} 
(1)$ the following statements are equivalent: \begin{enumerate}[$1)$]
\item $T$ is orthogonality preserving.
\item There exists a triple homomorphism $S: A \to
B^{**}$ satisfying that

$h^* S(a) =
S(a^*)^* h,$ $h S(a^*)^* = S(a) h^*,$ and \end{enumerate}

$$T(a) = h r(h)^* S(a)=
S(a) r(h)^* h,$$ for all $a\in A$.$\hfill \Box$
\end{theorem}

For an orthogonality preserving operator $T:A\to B$ we  call the triple homomorphism $S:A\to B^{**}$  given by Theorem \ref{t OPcstar}
the \emph{supporting triple homomorphism} of $T$. If $T$ is such that its supporting triple homomorphism is a TRO homomorphism, then we shall say that $T$ is a \emph{weighted TRO homomorphism}. Similarly, we define \emph{weighted TRO anti-homomorphisms}, \emph{weighted Jordan $^*$-homomorphisms}, etc.

We remark that from the proof of Theorem \ref{t OPcstar} (see the proof of Theorem 17 in \cite{BurFerGarMarPe}) the identities 
\begin{equation}\label{identities OPcstar}
 \begin{array}{ll}
         f(h^*) T(a) =
T(a^*)^* f(h), & f(h) T(a^*)^* = T(a) f(h^*)\\
        r(h)^* T(a) =
T(a^*)^* r(h) &\mbox{ and } r(h) T(a^*)^* = T(a) r(h)^*, \end{array}
\end{equation} hold for every $f\in C_0(Sp(h)), a \in A$.\smallskip

 The identities above are also satisfied when $T$ is replaced by its supporting triple homomorphism. We shall sometimes make use of these identities without explicit mention.\smallskip

 Let $A$ be a C$^*$-algebra. We denote by $Diag(a_1,a_2,\ldots,a_n)$ the element $(a_{i,j})$ in $M_n(A)$ such that $a_{i,j}=0$ if $i\neq j$ and $a_{i,i}=a_i.$ For $a\in A$ let us define $a_{_{(n)}}=Diag(a,a,\ldots,a)\in M_n(A).$ It is not hard to check that $(a_{_{(n)}})^{[2k-1]}=(a^{[2k-1]})_{_{(n)}}, k\in \NN.$ The usual arguments  show that $(a_{_{(n)}})^{\frac{1}{[2k-1]}}=(a^{\frac{1}{[2k-1]}})_{_{(n)}}$ $ ,k\in \NN.$ Taking weak$^*$-limits we have $r(a_{_{(n)}} )=r(a)_{_{(n)}}=Diag(r(a),r(a),\ldots,r(a)).$ If $r$ is a partial isometry, then straightforward computations allow to show that 
 $$r_{_{(n)}}r_{_{(n)}}^*M_n(B^{**}) r_{_{(n)}}^*r_{_{(n)}}=M_n(rr^*B^{**}r^*r).$$
        
 Let $T:A\to B, T=hr(h)^*S$ be an orthogonality preserving weighted TRO homomorphism.  Clearly we have $T_n=h_{_{(n)}} r(h)^*_{_{(n)}} S_n=S_n r(h)^*_{_{(n)}} h_{_{(n)}}.$ It is also not hard to check that $h_{_{(n)}}$ and $S_n$ satisfy the identities in Theorem \ref{t OPcstar} hence proving that $T_n$
 is orthogonality preserving, for every $n\in \NN.$ We next show that all c.o.p. operators between C$^*$-algebras arise this way.

  \begin{theorem}\label{t char cop} Let $T:A\to B$ be an orthogonality preserving operator between C$^*$-algebras. The following statements are equivalent:
  \begin{enumerate}[$1)$]

 \item $T$ is completely orthogonality preserving,

\item $T$ is 2-orthogonality preserving,

\item $T$ is a weighted TRO homomorphism,

\item $T$ preserves TRO-zero products,

\item $T(a)T(b)^*=0,$ whenever $ab^*=0.$

\end{enumerate}

  \end{theorem}

  \begin{proof}

  $1) \Rightarrow 2)$ It is obvious.

  $2)  \Rightarrow  3)$  Suppose that $T_{_2}$ is orthogonality preserving and let  $\widetilde{S}: M_2(A)\to M_2(B^{**})$ its supporting triple homomorphism. By the comments preceding the theorem we have $T_{_{2}}(1_{_{M_2(A^{**})}})=h_{_{(2)}},$ $r_2:=r(h_{_{(2)}})=r(h)_{_{(2)}}$ and if $x=\left[\begin{array}{ll}
         a & b\\
        c & d  \end{array} \right ]  \in M_2(A),$ then $T_{_2}(x)=h_{_{(2)}}r_{_{2}}^*\widetilde{S}(x). $ On the other hand we have
        \begin{equation}
        \begin{array}{c}\label{eq charCOP1}
        T_{_{2}} \left( \left[\begin{array}{ll}
         a & b\\
        c &d  \end{array} \right ] \right)=\left[\begin{array}{ll}
         hr(h)^*S(a) & hr(h)^*S(b)\\
        hr(h)^*S(c) &hr(h)^*S(d)  \end{array} \right ] =\\
         \left[\begin{array}{ll}
         h & 0\\
        0 &h  \end{array} \right ]  \left[\begin{array}{ll}
         r(h) & 0\\
        0 & r(h) \end{array} \right ]^*  \left[\begin{array}{ll}
         S(a) & S(b)\\
        S(c) & S(d)  \end{array} \right ]     =h_2 r_2^* S_2(x)   \end{array}
            \end{equation} 
 First observe that $\widetilde{S}(x)$ and $S_{_{2}}(x)$ both lie in $r_{_2}r_{_2}^*M_2(B^{**}) r_{_2}^*r_{_2}.$  By (\ref{eq charCOP1}) we have $ h_{_{(2)}}r_2^*(\widetilde{S}(x)-S_{_{2}}(x))=0.$ Similarly
        $ (\widetilde{S}(x)-S_{_{2}}(x))r_2^*h_{_{(2)}}=0.$  By triple functional calculus (see for instance $5)\Rightarrow 1)$ bellow) we have
        $$ r_2 r_2^*  (\widetilde{S}(x)-S_{_{2}}(x))=(\widetilde{S}(x)-S_{_{2}}(x)) r_2^* r_2=0. $$
          Since $ \widetilde{S}(x)-S_{_{2}}(x)$ lies in $r_2r_2^*M_2( B^{**})r_2^*r_2 $ we have that  $\widetilde{S}(x)=S_{_{2}}(x)$ for every $x\in M_2(A).$  
         We have proved that $S_{_{2}}$ is a triple homomorphism. By Proposition \ref{l 2-jordantriple} above, $S$ is
        a TRO homomorphism.

        $3) \Rightarrow  4)$  Let $T=hr^*S$ be a weighted TRO homomorphism. Recall that the ranges of both $S$ and $T$ are in  $ rr^*B^{**}r^*r.$ Let us fix $a,b,c\in A.$  By applying identities (\ref{identities OPcstar}) several times we obtain
         \begin{equation} \label{eq ThmCharCOP2}
        \begin{array}{c}
        T(a)T(b)^*T(c)=\\
        hr^*S(a)(hr^*S(b))^*S(c)r^*h=hr^*S(a)S(b)^*rh^*S(c)r^*h \\
        =hr^*S(a)S(b)^*S(c)r^*hr^*h=hr^*S(ab^*c)r^*hr^*h.  \end{array}
             \end{equation}
         As a consequence of (\ref{eq ThmCharCOP2}) if $ab^*c=0$  we have $T(a)T(b)^*T(c)=0,$ hence $T$ preserves TRO-zero products.\smallskip
         
         $4)\Rightarrow 5)$ The statement follows by Lemma \ref{lemma TRO-zero elements} since $T$ has the property $[a,b,b]=0 \Rightarrow [T(a),T(b),T(b)]=0$.\smallskip
         
         $5)\Rightarrow 1)$ Let us assume that $T$ preserves right orthogonality. Let us fix $a,b$ in $A$ with $ab^*=0.$ We have $0=T(a)T(b)^*=hr^*S(a)S(b)^*rh^*.$ From the last equalities we get $0=hh^*r^*S(a)S(b)^*rh^*=h^{[3]}r^*S(a)S(b)^*rh^*.$ Clearly, by induction, we have $0=h^{[2n-1]}r^*S(a)S(b)^*rh^*$ for all $n\in \NN.$ Hence $0=p(h)r^*S(a)S(b)^*rh^*$ for every odd polynomial. By the classical Stone-Weierstrass theorem we get $0=f(h)r^*S(a)S(b)^*rh^*$ for every $f\in C_0(Sp(h))$. In particular, the equality  $0=h^{[\frac{1}{3^n}]}r^*S(a)S(b)^*rh^*$
         holds. Taking w$^*$-limit  
         \begin{equation} \label{eq ThmCharCO3}
         0=rr^*S(a)S(b)^*rh^*.  
         \end{equation}
       
          Now, by using tripe functional calculus on $h^*$ we obtain from (\ref{eq ThmCharCO3})
          $$0=rr^*S(a)S(b)^*rr^*=rr^*S(a)(rr^*S(b))^*. $$
         Since $S(A)\subseteq rr^*B^{**}r^*r$ we finally have $S(a)S(b)^*=0.$ By Proposition \ref{char TRO-hom}, $S$ is a TRO homomorphism. Hence $T$ is c.o.p. (see the comments preceding the theorem).
         
\end{proof}

Let $A$ be an abelian C$^*$-algebra, $B$ a C$^*$-algebra and $T:A\to B$ be an orthogonality preserving operator. Let $S:A\to B^{**}$ be the supporting triple homomorphism of $T$. By Lemma \ref{l from abelian is TRO} we see that  $S$ is a TRO homomorphism and Theorem \ref{t char cop} then shows that $T$ is completely orthogonality preserving. Now if $A$ is a C$^*$-algebra such that for any C$^*$-algebra $B$ every orthogonality preserving operator $T:A\to B$ is completely orthogonality preserving, then $A$ is abelian. Indeed, in particular the mapping
$S:A\to A^{op}, \; a \mapsto a,$ is a Jordan $^*$-isomorphism which is c.o.p., by Proposition \ref{p triple2OP}, $S$ is a *-isomorphism. So $ab=ba, \forall a,b \in A,$ that is, $A$ is abelian.

\begin{corollary}
\label{l FROMabelian} Every orthogonality preserving operator from an abelian C$^*$-algebra is completely orthogonality preserving.  If $A$ is such that for any C$^*$-algebra $B$  every orthogonality preserving operator $T:A\to B$
is c.o.p., then $A$ is abelian. $\hfill \Box$
\end{corollary} 

Now we focus our attention on order zero maps and their relation with a problem posed in \cite{Gardner}. In the just quoted paper, zero product preserving operators are termed \emph{disjoint maps}. Theorem 1 in \cite{Gardner} states that a 2-positive, disjoint operator between C$^*$-algebras preserves absolute values. Under additional hypotheses 2-positive can be replaced by positive. The author leaves as an open problem whether the latter is true in general (see the comments after Corollary 7 in \cite{Gardner}). 
We give a positive answer to the problem posed by Gardner.

\begin{corollary}\label{c GardnerProblem} Let $T:A\to B$ be a positive  operator between C$^*$-algebras. The following statements are equivalent:
 \begin{enumerate}[$1)$]
 \item $T$ is a completely positive order zero map,

 \item $T$ preserves absolute values,

 \item $T$ preserves zero products,

 \item $T$ is 2-orthogonality preserving,

 \item $T$ is completely orthogonality preserving.
 \end{enumerate}

 \end{corollary}

 \begin{proof} $1)\Rightarrow 2)$ Follows from Theorem 2 in \cite{Gardner}.
$2)\Rightarrow 3)$ It follows  from the main theorem in \cite{Gardner} that $T$ is a weighted $^*$-homomorphism. Clearly, a weighted $^*$-homomorphism preserves zero products (take into account that in this case, $h$ commutes with the supporting $^*$-homomorphism of $T$).

 $3) \Rightarrow 4)$ This could be shown as a consequence of the non-unital version of Wolff's theorem contained in \cite{WinZac} (which, in turn, can also be obtained as a consequence of Wolff's theorem and Lemma 2.2 in \cite{ZeroWong}). We instead derive it from Theorem \ref{t OPcstar}. Since $T$ is positive and preserves zero-products then $T$ preserves orthogonality. By Theorem \ref{t OPcstar}, $T=hr(h)^* S.$ Since $h$ is positive then $p:=r(h)$ is a projection and $S(A)\subseteq pB^{**}p.$ As a consequence $pS=S$ and $S:A\to pB^{**}p$ is a Jordan $^*$-homomorphism. Since $pB^{**}p$ is a C$^*$-subalgebra of $B^{**}$ we easily realise that $S:A\to B^{**}$ is a Jordan $^*$-homomorphism. It also follows from Theorem \ref{t OPcstar} that $h$ commutes with $S.$ We claim that $S$  preserves zero products. Indeed, if $T(a)T(b)=0,$ then $h^2S(a)S(b)=0.$ As a consequence $h^n S(a)S(b)=0.$ By functional calculus (see for instance the proof of Lemma \ref{lemma TRO-zero elements}) it is not hard to see that this gives $p S(a)S(b)=0.$  Now if $ab^*=0$ then $S(a)S(b)^*=S(a)S(b^*)=0.$ By Proposition \ref{char TRO-hom}, $S$ is a TRO homomorphism. It follows from Theorem \ref{t char cop} that $T$ is $2$-orthogonality preserving.
 
 $4) \Rightarrow 5)$ It also follows from Theorem \ref{t char cop}.

$5) \Rightarrow 1)$ By Theorem \ref{t char cop}, $T=hr^*S,$ with $S$ being a TRO homomorphism. Following the same arguments given in  $3) \Rightarrow 4)$ we see that $S$ is a  $^*$-homomorphism, $r$ is a projection and $rS=S.$ For each natural $n$ and $x\in M_n(A)_{+}$ the elements $T_n(1_{M_n(A^{**})})=h_{_{(n)}}$ and $S_n(x)$ lie in $ M_n(B^{**})_{+}.$ Since  $h_{_{(n)}}$ and $S_n(x)$ commute, then $h_{_{(n)}} S_n(x)$ is also positive. Therefore $T$ is completely positive.
\end{proof}
  Let $T:A\to B$ be an orthogonality preserving operator with $T=hr^*S$. Then $T^{**}$ is an orthogonality preserving operator with supporting triple homomorphism $S^{**}.$  As a consequence $\|T \|=\|T^{**}\|=\| h\|.$  Further, if $T$ is c.o.p., then $\|T_n\|=\|h_n\|=\|h\|.$

 \begin{corollary}\label{c conOPisCC} Every contractive c.o.p. operator is completely contractive$\hfill \Box$\end{corollary} 
 
Let $S:M_m(\CC)\to M_m(\CC)$ be a TRO anti-homomorphism. Then $T=\frac{1}{m}S$
 is completely contractive and orthogonality preserving but it is not c.o.p.

\section{A generalisation of the order zero functional calculus}

Let $T:A\to B$ be an orthogonality preserving operator with $T=hr^*S$.  Let $h=T^{**}(1)$ and $f\in C_0(Sp(h)),$ where $Sp(h)$ stands for the triple spectrum of $h.$ We define the mapping
 $$f(T):A\to B^{**},\; f(T)(a):=f(h)r^* S(a).$$

 \begin{theorem}\label{t tripleFUNcalcForOP}
 $f(T)$ is an orthogonality preserving operator with $f(T)(A)\subseteq B$. Moreover, $f(T)$ is the unique orthogonality preserving operator with supporting triple homomorphism $S$ and such that $f(T)^{**}(1)=f(h).$
 
 If $T$ is c.o.p. then $f(T)$ is c.o.p.
 \end{theorem}
 \begin{proof}
We claim that the following identities
  \begin{equation}\label{identities TripleFUncCalc1}
 h^*T(a)^{^{[2n-1]}}=(T(a^*)^{^{[2n-1]}})^* h
      \mbox{ and } h (T(a^*)^{^{[2n-1]}})^*=T(a)^{^{[2n-1]}}h^*
\end{equation} hold for all $n\in \NN$ and $a\in A$.

 We shall only prove the first identity, the second one follows by similar arguments.\smallskip

 $h^*T(a)=T(a^*)^*h $ and $hT(a^*)^*=T(a)h^*$ by identities (\ref{identities OPcstar}).
Suppose that $ h^*T(a)^{^{[2k-1]}}=(T(a^*)^{^{[2k-1]}})^* h$ holds fo $k\leq n.$ Then $$ h^*T(a)^{^{[2n+1]}}=h^*T(a)^{^{[2n-1]}}T(a)^*T(a)=(T(a^*)^{^{[2n-1]}})^*h T(a)^*T(a)= $$ $$(T(a^*)^{^{[2n-1]}})^*T(a^*)h^*T(a)=(T(a^*)^{^{[2n-1]}})^*T(a^*)T(a^*)^*h=(T(a^*)^{^{[2n+1]}})^*h, $$ and the claim follows. Now we claim that \begin{equation}\label{identities TripleFUncCalc2}
 h^{^{[2n-1]}}r^* S(a^{^{[2n-1]}}) =T(a)^{^{[2n-1]}}, \mbox{ for all } n\in \NN,a\in A.
\end{equation} The case $n=1$ follows by Theorem \ref{t char cop}. Let us assume that $h^{^{[2k-1]}}r^* S(a^{^{[2k-1]}}) =T(a)^{^{[2k-1]}}$ for $k\leq n.$ Then
 $$h^{^{[2n+1]}}r^* S(a^{^{[2n-1]}})=h^{^{[2n+1]}}r^* S(a)^{^{[2n-1]}}=h^{^{[2n+1]}}r^* S(a)^{^{[2n-1]}}S(a)^*S(a)=$$
$$hh^*h^{^{[2n-1]}}r^* S(a)^{^{[2n-1]}}S(a)^*S(a)=hh^*T(a)^{^{[2n-1]}}S(a)^*S(a)= $$
 $$hh^*T(a)^{^{[2n-1]}}S(a)^*S(a)=h(T(a^*)^{^{[2n-1]}})^*hS(a)^*S(a)=T(a)^{^{[2n-1]}}h^*hS(a)^*S(a)= $$
 $$T(a)^{^{[2n-1]}}h^*hr^*rS(a)^*S(a)=T(a)^{^{[2n-1]}}T(a)^*T(a), $$ where the last equality is obtained after several applications of identities (\ref{identities OPcstar}). This proves (\ref{identities TripleFUncCalc2}).
 
Let us fix $a\in A,$ and $n$ natural. There exists an element $b\in A$ such that $b^{^{[2n-1]}}=a.$ By (\ref{identities TripleFUncCalc2}) we have $h^{^{[2n-1]}}r^* S(a)=h^{^{[2n-1]}}r^* S(b^{^{[2n-1]}}) =T(b)^{^{[2n-1]}},$ witnessing that $h^{^{[2n-1]}}r^* S(a)$ lies in $B.$
For any odd triple polynomial $p,$ the element $p(h)r^* S(a)$ is a linear combination of elements of the form $h^{^{[2n-1]}}r^* S(a) $ and hence $p(h)r^* S(a)$ lies in $B.$   Finally, by norm density of odd triple polynomials in $B^{**}_h$ we see that $f(h)r^* S(a)$ lies in $B,$ for every $f \in  C_0(Sp(h)).$ 
The rest of the proof follows easily from Theorem \ref{t OPcstar}

\end{proof}

 We call the mapping defined in the last theorem \emph{triple functional calculus on} $T.$ Let $T$ be an orthogonality preserving weighted TRO homomorphism. It is not hard to see that given $f_1,f_2,f_3\in C_0(Sp(h)),$ and $a,b,c\in A$ we have
 \begin{equation}\label{identites TROprodf(T)}
[f_1(T)(a),f_2(T)(b),f_3(T)(c)]=(f_1f_2^*f_3)(T)([a,b,c] ).
\end{equation}
 
 Let us fix $a\in A$ and let $(f_n)$  be a sequence in $ C_0(Sp(h))$ that converges in norm to $f\in C_0(Sp(h)).$ Then $$\|f_n(T)(a)-f(T)(a)\| \leq \| f_n(h)-f(h)\| \|a\|. $$
  This shows that the triple functional calculus on $T$ is continuous.
 
 Let now $T:A \to B$ be a symmetric orthogonality preserving operator between C$^*$-algebras. Then $h=T^{**}(1)$ is self-adjoint and $T=hJ,$ where $J:A\to B^{**}$ is a Jordan $^*$-homomorphism. In this case, we may define a functional calculus on $T$ by using the usual functional calculus for self-adjoint elements. Let C$^*(h)$ be the C$^{*}$-subalgebra of $B^{**}$ generated by $h,$ then C$^*(h)$ is $^*$-isomorphic to  $C_0(\Omega),$ with $\Omega$ a locally compact Hausdorff space.  Given $f$ in $C_0(\Omega)$ we define $f(T):=f(h)J.$ Similar (but simpler) computations to those in the proof of Theorem \ref{t tripleFUNcalcForOP} allow to prove:

 \begin{theorem}\label{t FUNcalcForSymmetricOP}
 $f(T)$ is an orthogonality preserving operator with $f(T)(A)\subseteq B$. Moreover, $f(T)$ is symmetric and is the unique orthogonality preserving operator with supporting Jordan $^*$-homomorphism $J$ and such that $f(T)^{**}(1)=f(h).$
 
  If $T$ is c.o.p. then $f(T)$ is c.o.p.$\hfill \Box$
 \end{theorem}

\begin{remark}\label{remark TROseminorm}
It is well known that the Gelfand-Naimark axiom is equivalent to the property $\|aa^*a\|=\|a\|^3.$ More generally, let $A$ be a $^*$-algebra and $\|.\|$ be a seminorm on $A$ which satisfies  $\|a b^*c\|\leq \|a\| \|b^*\| \|c\|$ and $\|aa^*a\|=\|a\|^3.$
 Since $\|a\|^3\leq \|a\|^2 \|a^*\|$ and $\|a^*\|^3\leq \|a^*\|^2 \|a\|$ then $\|a\|=0$ if and only if  $\|a^*\|=0.$ Thus it also follows that $\|a\|=\|a^*\|.$ It is not hard to check that
 \begin{equation}\label{eq TROseminorm}
     \|a\|^2=\sup \{ \|aa^*x\|: \|x\|\leq 1\}.
 \end{equation}
 
 holds. Let us fix $a \in A$ such that $\|a\|\neq 0.$ From (\ref{eq TROseminorm}) we have
  $\|a\|^6=\|aa^*a\|^2\leq \sup \{ \|aa^*\| \|aa^* \| \|aa^* x\| : \|x\|\leq 1\} \leq \| aa^*\|^2 \|a\|^2.$ Thus $\|a\|^2\leq \|aa^*\|.$ Finally we have $\| aa^*\|^3=\| aa^* aa^* aa^*\|\leq \|aa^*\| \|a\| \|a^*aa^*\|=\|aa^*\| \|a\|^4.$ This yields $\|aa^*\|\leq \|a\|^2.$ It follows that $\|.\|$ is a C$^*$-seminorm.
\end{remark}

  Let $E$ be a TRO and $h\in E$ with $\|h\|\leq 1.$ The TRO  $C_0((0,1])$ is generated by  the function $w(z)=z$ on $(0,1]$. The mapping $f\mapsto f_{|Sp(h)}$ is a TRO homomorphism between 
 $C_0((0,1])$ and $C_0(Sp(h)).$ Let $\psi^{-1}: C_0(Sp(h))\to E_h$ be the triple functional calculus on $h$. Then the
 mapping $\Phi_h: C_0((0,1]))\to E_h, f \mapsto \Phi_h(f)= \psi^{-1}(f_{|Sp(h)} ) $ is the unique (surjective) TRO homomorphism that extends the assignment $w\mapsto h.$ By an abuse of notation for $f\in C_0((0,1])$ we write $f(h)$ instead $\Phi_h(f).$ We shall use this convention without explicit mention during the proofs of the next results. We are now able to give a ternary version of \cite[Corollary 4.1]{WinZac}.

\begin{corollary}\label{cor BijectTROcop}
   \label{t PhiT}Let $T:A \to B$ be a contractive c.o.p. operator. Then there exists a TRO homomorphism  $\Phi_T:C_0((0,1])\otimes A\to B $ such that $\Phi_T(w\otimes  a)=T(a).$
 
 Conversely, any TRO homomorphism $\Phi: C_0((0,1])\otimes A\to B$ induces a contractive c.o.p. operator $T_{_{\Phi}}:A\to B$ via $T_{_{\Phi}}(a):=\Phi(w\otimes  a ).$
 
 Moreover,  the mappings $T\mapsto \Phi_{T}$ and $\Phi\mapsto T_{\Phi}$ are mutual inverses.\end{corollary}
\begin{proof} 
Since $T$ is  contractive then $\|h\|\leq 1 .$   For  $f\otimes a\in  C_0((0,1])\otimes A$ we set $\phi_T(f\otimes a):=f(T)(a)=f(h)r^*S(a).$ We can extend 
  $\phi_T$ to a linear mapping $\phi_T: C_0((0,1])\odot A\to B^{**}.$ By
   identity (\ref{identites TROprodf(T)}) and linearity we see that $\phi_T$ preserves TRO products. Actually, by the properties of the triple functional calculus on $T$ we have $\phi_T( C_0((0,1])\odot A)\subseteq B.$ Now let $x=\sum_i f_i \otimes a_i\in C_0((0,1])\odot A$ then 
$$\| \phi_T(\sum_i f_i \otimes a_i )\|\leq \sum_i \|f_i(h) r^*S( a_i)\|\leq \sum_i \|f_i\| \|a_i\|<\infty .$$   Thus the assignment $\| x\|_T:=\| \phi_T(x)\|$ defines a seminorm on  $ C_0((0,1])\odot A$ which enjoys the properties $\| xy^*z\|_T \leq \| x\|_T\| y^*\|_T \| z\|_T $ and $\| xx^* x\|_T=\| x\|^3_T.$ Thus $\| .\|_T$ is a C$^*$-seminorm by Remark \ref{remark TROseminorm}. Now recall that $C_0((0,1])$ is a nuclear C$^*$-algebra. We have $$\| \phi_T(x)\|=\| x\|_T\leq \| x\|_{max}=\|x\| $$
and therefore $\phi_T$ is contractive. Thus $\phi_T$ extends to a TRO homomorphism $\Phi_T:C_0((0,1])\otimes A\to B^{**}.$ Finally, $\Phi_T(C_0((0,1])\otimes A)\subseteq B$  and $\Phi_T(w\otimes a )=T(a)$ hold. 
 
  Conversely, if $\Phi: C_0((0,1])\otimes A\to B$ is a TRO homomorphism, then $T_{\Phi}:A\to B$ given by $T_{\Phi}(a)=\Phi(w \otimes a)$ is a c.o.p operator. Indeed,
  let $(a_{i,j}), (b_{i,j})$ in $M_2(A)$ with $(a_{i,j})\perp (b_{i,j}).$ Then we also have $w \otimes (a_{i,j})\perp w \otimes (b_{i,j})$ in $M_2(C_0((0,1])\otimes A).$ Since $\Phi$ is a TRO homomorphism then $\Phi_2:M_2(C_0((0,1])\otimes A)\to M_2(B^{**}) $ also is a TRO homomorphism. In particular, $\Phi_2$ is $2$-orthogonality preserving. As a consequence we have $T_2((a_{i,j}))\perp T_2((b_{i,j})).$ Hence $T$ is completely orthogonality preserving by Theorem \ref{t char cop}.\smallskip
  
  Finally, it is easy to see that the mappings $T\mapsto \Phi_{T}$ and $\Phi\mapsto T_{\Phi}$ are mutual inverses.
\end{proof}

 For $a_1,\ldots ,a_{2n+1}$ in a C$^*$-algebra  set $[a_1,\ldots,a_{2n+1}]:=[a_1,a_2,[a_3,\ldots, a_{2n+1}]],$ recursively.  Let $T:A\to B$ be a contractive c.o.p operator and $TRO(T(A))$ be the subTRO of $B$ generated by $T(A).$ 
  Then $TRO(T(A))$ is the norm closed linear span of all elements of the form $[T(a_1),T(a_2),\ldots,T(a_{2n+1})].$  
  
  \begin{corollary}\label{c ImagenOP}
  Let $T:A\to B$ be a contractive c.o.p. operator. Then
  $$TRO(T(A))= \Phi_T(C_0((0,1])\otimes A). $$
  \end{corollary}
  \begin{proof}
  Similar computations to those in the proof of Theorem \ref{t tripleFUNcalcForOP}  show that
  \begin{equation}\label{eq TROTA1}
      [T(a_1),\ldots,T(a_{2n+1})]=w^{[2n+1]}r^*S([a_1,\ldots,a_{2n+1}])\end{equation} (recall the convention established in the comments previous to Corollary \ref{cor BijectTROcop}).
  By (\ref{eq TROTA1}) every element of the form $p(T)r^*S(a),$ where $p$ is an odd polynomial, lies in $TRO(T(A)).$ Now let us fix $ f\in C_0((0,1]).$ There exists a sequence $(p_n)$ of odd polynomials such that $(p_n)\to f$ in norm. By continuity of the triple functional calculus on $T$ we have 
  $(p_n(T)(a))  \to f(T)(a),$ for every $a\in A.$ Since $TRO(T(A))$ is norm-closed, this shows that every element of the form $f(T)(a)$ also lies in $TRO(T(A)).$  Let us write $$C=\overline{span\{f(T)(a): f\in C_0((0,1]),a\in A \} }^{\|. \|}$$  It follows that 
  $C\subseteq TRO(T(A)).$ The other inclusion is straightforward. Hence we have $TRO(T(A))=C.$
   
  Now let $\Phi_T:C_0((0,1])\otimes A\to B$ be the TRO homomorphism induced by $T$ (see Theorem \ref{t PhiT}). First recall that, since $\Phi_T$ is a triple homomorphism, then $\Phi_T$ has closed range (compare with \cite[Corollary 7]{KapsTriples}) and hence $\Phi_T(C_0((0,1])\otimes A)$ is a TRO.
  Let $\sum_{n=1}^{\infty} f_n\otimes a_n$ be an element in $ C_0((0,1])\otimes A$. Then
  \begin{equation}\label{eq TROTA3}
   \Phi_T(\sum_{n=1}^{\infty} f_n\otimes a_n)=\sum_{n=1}^{\infty} f_{n}(T)(a_n)   
  \end{equation}
   which lies in $TRO(T(A)).$ For the other inclusion observe that since
    $TRO(T(A))=C$ then any element in $TRO(T(A))$ is of the form $\sum_{n=1}^{\infty} f_{n}(T)(a_n)$ which, by (\ref{eq TROTA3}), lies in $\Phi_T(C_0((0,1])\otimes A).$
   \end{proof}

  Let $T:A\to B$ be a positive c.o.p. operator (equivalently, a c.p. order zero map by Corollary \ref{c GardnerProblem}). In this case $T$ is a weighted $^*$-homomorphism by Corollary \ref{c GardnerProblem} and induces a $^*$-homomorphism $\Psi_T:C_0((0,1])\otimes A\to B$  such that $\Psi_T(w\otimes a )=T(a) $ (see also \cite[Corollary 4.1]{WinZac}). Similar computations to those performed in the proof of Corollary \ref{c ImagenOP} allow to prove:
  
  \begin{corollary} If $T:A \to B$ is a c.o.p. positive operator then
  $C^*(T(A))= \Psi_T(C_0((0,1])\otimes A ).$ $\hfill \Box$
  \end{corollary}
 
\renewcommand{\abstractname}{Acknowledgements}
\begin{abstract}
 The author would like to thank the anonymous reviewers for their careful reading and for several comments that allowed to improve the paper.
\end{abstract}


\begin{thebibliography}{0}

\bibitem{AlaminosZeroP} J. Alaminos, M. Bre\v sar, J. Extremera and A. Villena, Maps preserving zero products, \emph{Studia Math.} \textbf{193} no. 3 (2009) 131-159.

\bibitem{Arendt} W. Arendt, Spectral properties of Lamperti operators,
\emph{Indiana Univ. Math. J.}, \textbf{32} no. 2 (1983) 199-215.

\bibitem{BarDangHorn} T. J. Barton, T. Dang and G. Horn, Normal representations of Banach Jordan triple systems, \emph{Proc. Amer. Math. Soc.}, \textbf{102} no. 3 (1988) 551-555.


\bibitem{BeckNarici}  E. Beckenstein, L. Narici, Automatic continuity of certain linear isomorphisms, \emph{Acad. Roy. Belg. Bull. Cl. Sci.} \textbf{73} no. 5 (1987) 191-200.



\bibitem{BurFerGarMarPe} M. Burgos, F.J. Fern{\' a}ndez-Polo,
J.J. Garc{\'e}s, J. Mart{\'i}nez Moreno, A.M. Peralta,
Orthogonality preservers in C*-algebras, JB*-algebras and
JB*-triples, \emph{J. Math. Anal. Appl.}, \textbf{348} (2008) 220-233.

\bibitem{BurFerGarPe} M. Burgos, F.J. Fern{\' a}ndez-Polo, J.J. Garc{\'e}s, A.M. Peralta, Orthogonality preservers Revisited,
\emph{Asian-European J. Math.}, \textbf{2} no. 3 (2009) 387-405.

\bibitem{bussim} A. Buss and A. Sims, Opposite algebras of groupoid C$^*$-algebras, preprint.


\bibitem{CheKeLeeWong} M. A. Chebotar, W.F. Ke, P.H. Lee, and N.C. Wong, Mappings preserving zero products, \emph{Studia Math.}, \textbf{155} no. 1 (2003) 77-94.

\bibitem{connes}  A. Connes, A factor not anti-isomorphic to itself, \emph{Ann. Math.},  \textbf{101} no. 2 (1975)  536-554.

\bibitem{EdRut} C.M. Edwards, G.T. Ruttimann, Compact tripotents in bi-dual JB$^*$-triples, \emph{Math. Proc. Cambridge Philos. Soc.}, \textbf{120} no. 1 (1996)  155–173. 

\bibitem{KapsTriples} F.J. Fern{\' a}ndez-Polo, J.J. Garc{\'e}s, A.M. Peralta, A Kaplansky Theorem for JB$^*$-triples, \emph{Proc. Amer. Math. Soc.}, \textbf{140} no. 9 (2012) 3179-3191.

\bibitem{franzoni} T. Franzoni, The group of holomorphic automorphisms in certain J$^*$-algebras, \emph{Ann. Mat. Pura Appl.},  \textbf{CXXVII} no. IV (1981) 51-56.
 

\bibitem{Gardner0} L.T. Gardner, A dilation theorem for $|.|$-preserving maps of C$^*$-algebras, \emph{Proc. Amer. Math. Soc.}, \textbf{73} (1979) 461-465.

\bibitem{Gardner} L.T. Gardner, Linear maps of C$^*$-algebras preserving the absolute value, \emph{Proc. Amer. Math. Soc.}, \textbf{79} no. 2 (1979) 271-278.

\bibitem{HAM} M. Hamana, Triple envelopes and Silov boundaries of operator spaces, \emph{Math. J. Toyama Univ.}, \textbf{22} (1999) 77–93.

\bibitem{HaOlSt}H. Hanche-Olsen, E. Stormer, Jordan operator algebras. Pitman, London (1984).

\bibitem{Harris} L.A. Harris, A generalization of C$^*$-algebras, \emph{Proc. Lond. Math. Soc.}, \textbf{42} no. 3 (1981) 331-361.

\bibitem{Horn} G. Horn,
Characterization of the predual and ideal structure of a JBW$^*$-triple,
\emph{Math. Scand.}, \textbf{61} no. 1 (1987) 117-133.

\bibitem{Jarosz}  K. Jarosz,  Automatic  continuity  of  separating  linear  isomorphisms, \emph{Canad. Math. Bull.}, \textbf{33} no. 2 (1990) 139-144 .

\bibitem{JeangWong} J.S. Jeang, N.C. Wong, Weighted composition operators of
$C_0(X)$'s, \emph{J. Math. Anal. Appl.} \textbf{201},  981-993 (1996).

\bibitem{Kads} R.V. Kadison, Isometries of Operator algebras, \emph{ Ann. Math.}, \textbf{54} no. 2 (1951) 325-338.

\bibitem{Kaup} W. Kaup, A Riemann mapping theorem for bounded symmetric domains in complex Banach spaces, \emph{Math. Z.}, \textbf{183} (1983) 503-529. 

\bibitem{LauWong} A.T.M. Lau and N.C. Wong, Orthogonality and disjointness preserving linear maps between Fourier and Fourier-Stieltjes algebras of locally compact groups, \emph{J. Funct. Anal.}, \textbf{265} (2013) 562-593.

\bibitem{JuChunChinWong}J.H. Liu, C.Y. Chou, C.J. Liao and N.-C. Wong, Linear disjointness preservers of operator algebras and related structures. \emph{Acta Sci. Math.}, \textbf{84} no. 1-2 (2018) 277-307.

\bibitem{PeraltaLocalRep} A. M. Peralta, A note on 2-local representations, \emph{Operators ans Matrices.}, \textbf{9} no. 2 (2015) 343-358.
 
\bibitem{sakai} S. Sakai, C$^*$-algebras and
W$^*$-algebras. Springer Verlag. Berlin (1971).

\bibitem{takesaki} M, Takesaki, Theory of Operators algebras I, Springer Verlag, Berlin, 2002.

\bibitem{tomiyama} J. Tomiyama, On the transpose map of matrix algebras, \emph{Proc. Amer. Math. Soc.}, \textbf{88} no. 4 (1983) 635-638.

\bibitem{WinterCovering} W. Winter, Covering dimension for nuclear $C^*$-algebras, \emph{J. Fun. Anal.}, \textbf{199}, no. 2 (2003) 535-556.

\bibitem{WinZac} W. Winter and  J. Zacharias. Completely positive maps of order zero. \emph{M\"unster J. of Math.}, \textbf{2} (2009) 311-324.

\bibitem{Wolff} M.Wolff, Disjointness preserving operators on C$^*$-algebras, \emph{Arch. Math.} \textbf{62} (1994) 248-253.

\bibitem{TripleWong} N.C. Wong, Triple homomorphisms of C$^*$-algebras, \emph{Southeast Asian Bulletin of Math.}, \textbf{29} no. 2 (2005) 401-407.

\bibitem{ZeroWong} N.C. Wong, Zero product preservers of C$^*$-algebras, \emph{Contemporary Math.}, \textbf{475} (2007), 377-380.

\bibitem{Zetll} H. Zettl, A characterization of ternary rings of operators, \emph{Adv. in Math.} \textbf{48} (1983) 117-143.

\end{thebibliography}
 \end{document}